# Accelerating the Serviceability-Based Design of Reinforced Concrete Rail Bridges under Geometric Uncertainties induced by unforeseen events: A Surrogate Modeling approach

Mouhammed ACHHAB, Pierre JEHEL, Fabrice GATUINGT

Université Paris-Saclay, CentraleSupélec, ENS Paris-Saclay, CNRS, LMPS – Laboratoire de Mécanique Paris-Saclay, 91190, Gif-sur-Yvette, France,
Mouhammed Achhab (mouhammed.achhab@centralesupelec.fr)

## Abstract

Reinforced concrete rail bridges are essential components of railway infrastructure, where reliability, durability, and adaptability are key design priorities. However, the design process is often complicated by uncertainties stemming from unforeseen construction constraints, such as the need to reposition piers or alter geometric characteristics. These design adaptations can lead to repeated redesigns, added costs, and project delays if not anticipated in the early design stages, as well as significant computational overhead when using traditional finite element (FE) simulations. To address this and anticipate such unexpected events, this study adopts surrogate modeling as an efficient probabilistic design approach. This methodology integrates key geometric parameters as random variables, capturing the uncertainties that may arise during the design and construction phases and propagating them on the bridge's performance functions. By doing so, we aim to enable the efficient exploration of a large number of design scenarios with minimal reliance on time-consuming finite element (FE) simulations, represent the performance functions of a reinforced concrete bridge as a function of our variable design parameters, and classify the overall design scenarios into failure and safe scenarios In this study, a four-span reinforced concrete bridge deck is modeled using a multi-fiber finite element approach in Cast3M software. This FE model is used to generate the required design of experiments to train the surrogate models. Within this framework, a comparative performance assessment is conducted to evaluate the performance of the Kriging surrogate against alternative methods, including polynomial chaos expansion (implemented in UQLab) and support vector regression (SVR). This methodology supports early-stage uncertainty-informed design, enhancing the robustness and adaptability of reinforced concrete rail bridges in the face of practical constraints and changing site conditions.

## 1. Introduction

The growing demand for rail infrastructure reflects a global push for faster, more connected, and efficient transportation systems, particularly in response to urban expansion and population density. Reinforced concrete rail bridges have become a common and reliable structural choice, valued for their adaptability, durability, and cost-effectiveness. Yet, their design and construction remain complex due to large scale, variable operating conditions, and the nonlinear behavior of concrete.

Traditionally, engineers rely on finite element (FE) models combined with standardized design practices, such as those defined in the Eurocodes (European Committee for Standardization (CEN), 2004), to ensure that both safety and serviceability requirements are satisfied. These requirements are framed around the limit state design methodology, which distinguishes between ultimate limit states (ULS), related to structural failure, and serviceability limit states (SLS), which ensure usability and comfort. Examples of ULS include reaching flexural or shear strength, whereas SLS are typically defined by constraints on deflection, cracking, or vibration.

Design codes account for uncertainties through partial safety factors, which correct for variability in loads, material properties, and construction tolerances. Yet, real-world uncertainties often arise during the design and construction phases that lie outside the scope of predefined code factors. For example, unforeseen underground obstructions or the need to relocate structural elements during construction may force changes that invalidate initial assumptions, triggering redesign cycles, delays, and increased costs.

This study addresses such geometrical uncertainties by proposing a probabilistic framework aimed at improving the robustness and adaptability of reinforced concrete bridge designs. Specifically, we consider uncertainties stemming from unexpected events that affect the geometry of the bridge, such as pier positioning, which may vary due to site-specific constraints. To support efficient exploration of this uncertain design space, we employ a surrogate modeling approach that reduces the computational burden typically associated with high-fidelity FE simulations
.
Our methodology involves constructing a predictive model based on a set of multifiber beam (Combescure 2007) finite element simulations of a continuous reinforced concrete bridge deck. This surrogate model approximates the behavior of the structure across varying geometrical configurations, allowing for fast evaluation of a limit state function under uncertainty. Unlike conventional analyses focused on deflection limits, we target, in the present case, a different serviceability criterion: the maximum tensile stress in steel reinforcement. This quantity is directly related to crack width control, a critical consideration for long-term durability and maintenance. The uncertain parameters are related to the locations of the three internal piers of our bridge deck and the thickness of the bridge deck. This approach facilitates the early-stage integration of probabilistic analysis into the design process, aiming to reduce redesign iterations and improve resilience to unforeseen construction challenges.

The Kriging (Sacks et al. 1989) or Gaussian Process (Rasmussen and Williams 2008) method is selected as a primary surrogate model due to its accuracy and flexibility, and its proven performance in uncertainty propagation, design optimization (Dubourg, Sudret, and Bourinet 2011), and reliability analysis (Lelièvre et al. 2018). Within this framework, a comparative performance assessment is also conducted to compare the performance of Kriging (Lataniotis, Marelli, and Sudret 2018) with different types of regression functions, with other methods, including polynomial chaos expansion (PCE) (Xiu and Karniadakis 2002) and Support Vector Regression (SVR) (Smola and Schölkopf 2004), across various dataset sizes and evaluation metrics.

Section II describes the finite element model of the reinforced concrete bridge used as the basis for surrogate model development. Section III outlines the surrogate modeling approach adopted in this study and its numerical implementation. In Section IV, the Results are presented and analyzed. Finally, Section V summarizes the main conclusions and highlights directions for future research.

## 2. Parametric Finite element model of the bridge deck

The finite element model is developed using Cast3M, a finite element software well-suited for structural and nonlinear analyses. The structure under study is a four-span reinforced concrete bridge deck with an overall length of 40 m and a constant width of 8 m. This bridge deck is supported on five piers, which prevent both vertical and horizontal displacements. The two end supports are positioned at fixed locations: 0 m and 40 m along the bridge length. The applied loading follows the LM-71 model from Eurocode 1, doubled to represent two loaded lanes in accordance with multi-lane traffic design rules. All geometrical parameters of the bridge deck and the loading conditions are illustrated in **Figure 1**.

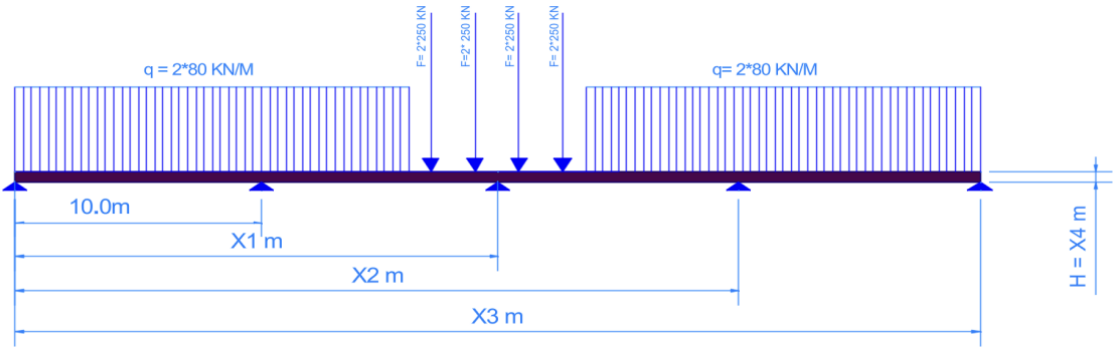

**Figure 1.a:** Longitudinal view of the bridge deck

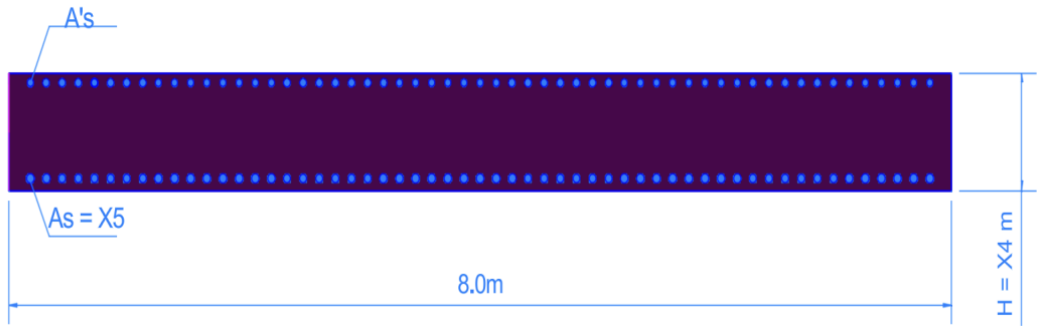

**Figure 1.b:** Transversal section of the bridge deck

In this study, a selected subset of geometrical design parameters is considered as sources of uncertainty. Specifically, the bridge deck thickness and the longitudinal positions of the three intermediate piers are modeled as independent random variables. While the steel reinforcement ratio is not directly treated as a random input, it is implicitly affected by the variability in deck thickness, given its definition as a fixed percentage of the concrete cross-sectional area. The uncertainty is described by the random vector X = [$X_1$, $X_2$, $X_3$, $X_4$], where each component represents a distinct geometrical quantity. $X_1$ corresponds to the position of the second pier, $X_2$ captures the variability in the position of the third (central) pier, $X_3$ represents the position of the fourth pier, and $X_4$ models the deck thickness. The geometrical interpretation of each of these parameters is illustrated in **Figure 1**, while the probabilistic characterization of their distributions (type, range, and parameters) is provided in **Table 1**.

| Type | Geometric | Geometric | Geometric | Geometric |
|---|---|---|---|---|
| Design Parameter | Position of 2nd pier | Position of 3rd pier | Position of 4th pier | Bridge Deck thickness |
| Symbol | $X_1$ | $X_2$ | $X_3$ | H |
| Random Distribution | $X_1 \sim U(5m,15m)$ | $X_2 \sim U(15m,25m)$ | $X_3 \sim U(25m,35m)$ | $X_4 \sim U(30cm,60cm)$ |

**Table 1:** Random distributions of geometrical design parameters

A multifiber beam formulation is employed in this study as a compromise between classical beam elements and fully 3D volumetric models, offering a balance between computational efficiency and the ability to capture detailed structural behavior. In the longitudinal direction, the bridge deck is discretized using Timoshenko beam elements (Öchsner 2021), allowing for the representation of both bending and shear deformations. Each Gauss integration point along these longitudinal elements corresponds to a transverse section, which is then discretized into concrete and steel fibers. This fiber-level discretization, applied at the Gauss points, is essential for assigning appropriate material-specific constitutive laws to each fiber. In particular, the QUAS model with four integration points is used for concrete fibers to capture nonlinear and damage behaviors, while the POJS model with a single integration point is applied to the steel fibers to represent their elasto-plastic behavior. To ensure high-fidelity results, a fine mesh is used in the simulations, with a density parameter of 0.5 in Cast3M, providing a detailed resolution of both global structural response and local stress-strain evolution across sections. For the constitutive modeling of materials, a nonlinear approach is adopted. The behavior of concrete is simulated using the De Laborderie model, implemented in Cast3M under the Unilateral formulation. This model captures the uniaxial behavior and damage mechanisms inherent in concrete under tension and compression. Steel reinforcement is modeled with plastic hardening, thereby allowing the simulation of strain-hardening effects, essential for realistic structural analysis under high loading scenarios.

## 3. Surrogate modeling

The serviceability limit state related to maximum stress in the steel rebars is defined by the function $g(x) = L - q(x)$, where L denotes the maximum allowable stress in steel rebars threshold, and $q(x)$ represents the maximum stress in steel rebars obtained from the finite element model for a given set of design parameters $x$. The bridge is considered to remain in serviceable condition as long as $g(x) > 0$. To develop the surrogate model based on FEM simulations we need the necessary input-output data set, so n samples $S = \{x^{(1)}, ..., x^{(n)}\}$ should be selected using LHS and the output of each input point $g^{(i)}(x) = g(x^{(i)})$ should be calculated to create the output set $\mathcal{G} = \{g(x^{(1)}) ... g(x^{(n)})\}$. To this end, as shown in **Figure 2**, a computational workflow was established, integrating MATLAB and Cast3M. Latin Hypercube Sampling (LHS) is utilized within MATLAB to generate input samples for the design variables $x$. For each sampled input, MATLAB automatically triggers Cast3M to perform the finite element analysis and compute the corresponding maximum stress in steel rebars. The resulting outputs, along with their associated input values, are collected and returned to MATLAB, where post-processing, data aggregation, and surrogate modeling are carried out.

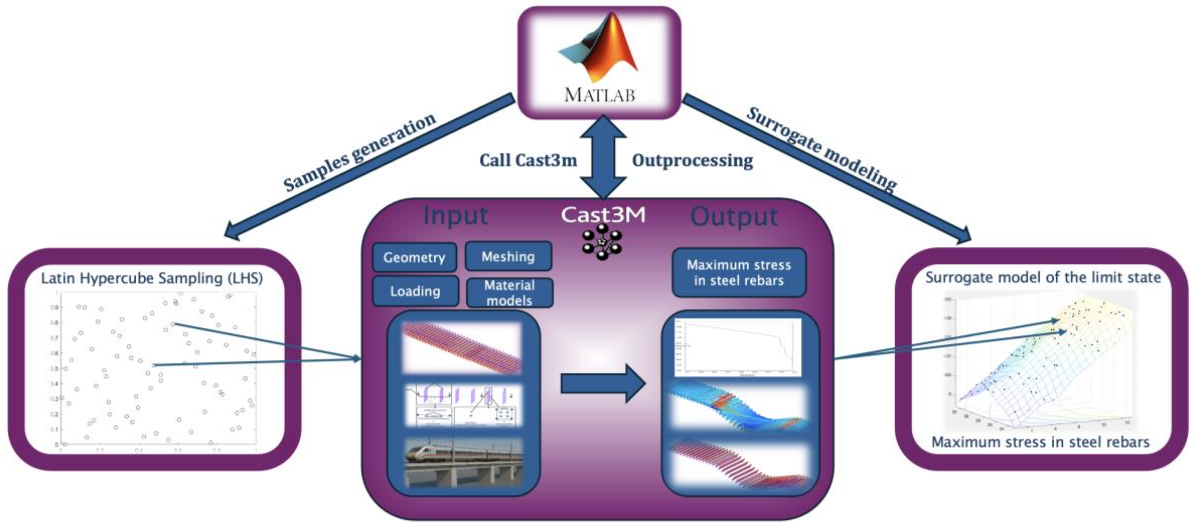

**Figure 2:** Cast3m/MATLAB computational pipeline

To implement our methodology and develop the surrogate models, we proceeded as follows. A set of 50 simulation results was set aside as a validation set, while multiple training datasets of increasing size (50, 100, ..., 300 simulations) were gathered in parallel. For each training set, several surrogate models were developed using three different techniques: Kriging, Polynomial Chaos Expansion (PCE), and Support Vector Regression (SVR). In the case of Kriging, five surrogate models were constructed: one using ordinary Kriging (with a constant trend), and four others incorporating polynomial trends of degree 1 to 4. For PCE, a single model was built, with the polynomial degree automatically selected by the algorithm. One SVR model was also developed for each training dataset. The performance of all surrogate models was evaluated using multiple error metrics, including the Mean Absolute Error (MAE), Maximum Absolute Error (MaxAE), Root Mean Square Error (RMSE), and the coefficient of determination ($R^2$). The primary objective of this study is to assess the performance of the Kriging-based surrogate models and compare their accuracy against those obtained using PCE and SVR.

## 4. Results

**Figure 3** presents the evolution of four performance metrics—MAE, MaxAE, RMSE, and $R^2$—used to assess the accuracy of several surrogate models with varying training set sizes. The surrogates compared include Ordinary Kriging, Universal Kriging with linear, quadratic, third-order, and fourth-order polynomial regression functions, as well as Polynomial Chaos Expansion (PCE) and Support Vector Regression (SVR). MAE evaluates the average error magnitude, MaxAE highlights the worst-case (local) error, RMSE emphasizes larger deviations due to its squared formulation, and $R^2$ measures the global predictive capability, ranging from 0 to 1, with higher values indicating a better fit. In our study, $R^2$ is primarily used to quantify the overall performance of each surrogate, complemented by the other metrics to assess both global and local accuracy, as well as error trends across methods.

Based on the results, Universal Kriging with a third-order polynomial regression function consistently outperformed all other models. It achieved an $R^2$ close to 0.99 with 150 training samples, and exceeded 0.99 with 200 samples. It also maintained the lowest RMSE and MAE across all training set sizes. With respect to MaxAE, which captures local error peaks, the third-order Kriging model also yielded the lowest or near-lowest values, demonstrating robust local accuracy. Kriging performance decreased when using a regression order higher or lower than 3, or when reverting to Ordinary Kriging (i.e., no regression term), indicating that the third-order polynomial trend offers an optimal level of accuracy.

Most surrogates reached high accuracy with 300 training samples, achieving $R^2$ values in the 0.98–0.99 range. The exception was PCE, which initially showed good performance (especially with as few as 50 samples), but plateaued and then declined slightly beyond 200 samples, ending with an $R^2$ of approximately 0.96. Nevertheless, PCE outperformed Kriging models with linear, fourth-order regression, and Ordinary Kriging for smaller datasets (≤200 samples).

SVR demonstrated fluctuating performance: it improved up to 150 training samples, then declined at 200 and 250, before regaining high accuracy at 300 samples. This indicates some instability in SVR's response to increasing data, despite reaching good performance levels at larger sample sizes.

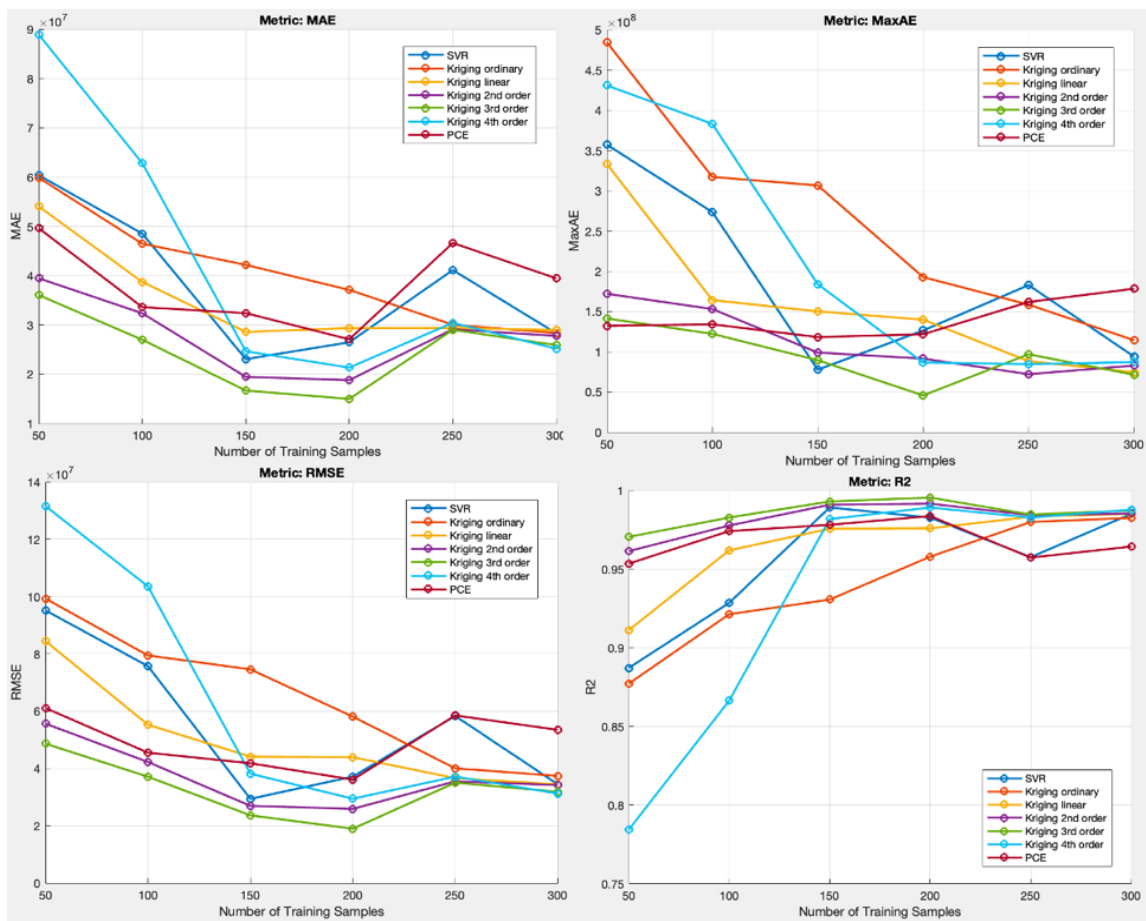

**Figure 3:** Performance of the surrogates

## 5. Conclusion

The results demonstrate that surrogate modeling is highly effective for our problem, enabling accurate prediction of the maximum stress in steel reinforcement with significantly reduced computational cost. Among all tested models, Kriging-based surrogates delivered the most consistent and accurate performance, especially the Universal Kriging with a third-order polynomial trend, which outperformed others across all evaluation metrics.

Beyond accuracy, Kriging offers strong interpretability, as it provides both a prediction and a measure of uncertainty at each input point, allowing better insight into model confidence and guiding further simulation efforts. The comparison also highlights that the choice of regression trend in Kriging plays a crucial role in model accuracy.

Overall, Kriging is shown to be the most reliable and interpretable surrogate model for this application, with trend selection being a key aspect of its successful implementation.

Future work will explore multi-fidelity and active learning approaches to reduce computational cost and improve surrogate accuracy near the limit-state surface. Additionally, extending the analysis to multiple limit states, such as deflection alongside reinforcement stress, could be explored within an active learning approach.

## Acknowledgment

This research is funded by the MINERVE project. The MINERVE project has been financed by the French government within the framework of France 2030.